\documentclass[a4paper]{article}
\usepackage{amscd,amsmath,amsthm,amssymb}
\begin{document}

\newtheorem{theorem}{Theorem}[section]
\newtheorem{lem}[theorem]{Lemma}
\newtheorem{prop}[theorem]{Proposition}
\newtheorem{cor}[theorem]{Corollary}

\theoremstyle{definition}
\newtheorem{definition}[theorem]{Definition}
\newtheorem{example}[theorem]{Example}

\theoremstyle{remark}
\newtheorem*{remark}{{\sc Remark}}

\newcommand{\debar}{{\overline{\partial}}}
\newcommand{\mapor}[1]{{\stackrel{#1}{\longrightarrow}}}
\newcommand{\mapver}[1]{\Big\downarrow\vcenter{\rlap{$\scriptstyle#1$}}}
\newcommand{\sA}{{\mathcal A}}
\newcommand{\Oh}{{\mathcal O}}
\newcommand{\sH}{{\mathcal H}}
\newcommand{\K}{\mathbb{K}}
\newcommand{\R}{\mathbb{R}}
\newcommand{\C}{\mathbb{C}}
\newcommand{\Z}{\mathbb{Z}}
\newcommand{\N}{\mathbb{N}}
\newcommand{\Coker}{\operatorname{Coker}}
\newcommand{\Hom}{\operatorname{Hom}}
\newcommand{\Image}{\operatorname{Im}}

\title{Deformation theory via differential graded Lie algebras}
\author{Marco Manetti,
Universit\`a di Roma ``La Sapienza'', Italy}
\date{11 July 2005}

\maketitle

\def\abstractname{Read me first}

\begin{abstract}
This is an expository paper written in 1999 and published in
\emph{Seminari di Geometria Algebrica 1998-1999,}
Scuola Normale Superiore (1999).\\
Six years later some arguments used here appear quite naive and
obsolete but, in view of the several citations that this paper has
obtained in the meantime,  I preferred don't change the
mathematical contents and to fix only some typos and minor
mistakes.\\
For a better and more recent  treatment of these topics we refer
to: M. Manetti, \emph{Lectures on deformations on complex
manifolds,} Rendiconti di Matematica \textbf{24} (2004) 1-183.
\end{abstract}

\section*{Introduction}

This paper concerns the basic philosophy that, over a field of
characteristic 0, every deformation problem is governed by a differential
graded Lie algebra (DGLA) via solutions of Maurer-Cartan equation modulo
gauge action.\\
The classical approach (Grothendieck-Mumford-Schlessinger) to
infinitesimal deformation theory is described by the procedure (see e.g
\cite{Ar})
\[\hbox{Deformation problem}\quad\leadsto\quad\hbox{Deformation
functor}\]%
The above picture is rather easy and suffices for many
applications; it is however clear that in this way we forget
information which can be useful. The other classical approach,
which consider categories fibred in groupoids instead of
deformation functors,
is not much better.\\
A possible and useful way to preserve information is to consider a
factorization
\[\hbox{Deformation problem}\quad\leadsto\quad
DGLA\quad\leadsto\quad \hbox{Deformation functor}\]%
where by $DGLA$ we mean a differential graded Lie Algebra
depending from the data of the deformation problem and the arrow
$DGLA\;\leadsto\; \hbox{Deformation functor}$, is a well defined
and functorial procedure explained in Section \ref{sec.3}.
Moreover we prove that every quasiisomorphism of differential
graded Lie algebras induces an
isomorphism of deformation functors.\\
Given a deformation problem, in general it is not an easy task to
find a factorization as above; some general technics of this
``art'' (see \cite[p. 5]{K2}), will be discussed elsewhere. Here
we only point out that in general the correct DGLA is only defined
up to quasiisomorphism and then this note represents the
necessary background for the whole theory.\\
Although the interpretation of deformation problems in terms of solutions
of Maurer-Cartan equation is  very useful on its own,
in many situation it is unavoidable to recognize that
the category of DGLA is too rigid for a ``good'' theory.
The appropriate way of extending this category will be the
introduction of homotopy Lie algebras and $L_{\infty}$-algebras; these new
objects will be described in next lectures.\\
The results of this paper are, more or less, known to
experts; if some originality is present in this notes then it is only
contained in the proofs.

\section{Differential graded Lie algebras and the Maurer-Cartan equation}
\label{sec.1}

Unless otherwise specified we shall follow the notation of \cite{K1}.\\
$\K$ is a field of characteristic 0.\\
$Graded^{\K}$ is the category of $\Z$-graded vector space over
$\K$.\\
$[n]\colon Graded^{\K}\to Graded^{\K}$ is the shift operator,
$V[n]^i=V^{i+n}$; for example $\K[1]^{-1}=\K$, $\K[1]^i=0$ for
$i\not=-1$.\\
If $V$ is a $\Z$-graded vector space and
$v\in V$ is a homogeneous element we
write $\overline{v}\in\Z/2\Z$ for the class modulo 2 of the
degree of $v$.\\

\begin{definition}\label{1.1} A \emph{Differential graded Lie algebra} (DGLA for short)
$(L,[,],d)$ is the data of a
$\Z$-graded vector space
$L=\oplus_{i\in \Z}L^i$ together a bilinear bracket
$[,]\colon L\times L\to L$ and a linear map $d\colon L\to L$ satisfying
the following condition:\begin{enumerate}

\item $[,]$ is homogeneous skewsymmetric;
this means $[L^i,L^j]\subset L^{i+j}$ and
$[a,b]+(-1)^{\overline{a}\overline{b}}[b,a]=0$ for every $a,b$ homogeneous.

\item  Every $a,b,c$ homogeneous satisfy the Jacobi identity
\[ [a,[b,c]]=[[a,b],c]+(-1)^{\overline{a}\overline{b}}[b,[a,c]].\]

\item  $d(L^i)\subset L^{i+1}$, $d\circ d=0$ and
$d[a,b]=[da,b]+(-1)^{\overline{a}}[a,db]$. The map $d$ is called the
differential of $L$.
\end{enumerate}
\end{definition}

Note that $L^0$ and $L^{even}=\oplus L^{2n}$ are Lie algebras in the
usual sense.\\

\begin{definition}\label{1.2}
A linear map $f\colon L\to L$ is called a {\it derivation of
degree $n$} if $f(L^i)\subset L^{i+n}$ and satisfies the graded
Leibniz rule $f([a,b])=[f(a),b]+(-1)^{n\overline{a}}[a,f(b)]$.\end{definition}

We note that if $a\in L^i$ then $ad(a)\colon L\to L$, $ad(a)(b)=[a,b]$,
is a derivation of degree $i$ and $d$ is a derivation of degree 1.\\

By following the standard notation we denote by $Z^i(L)=\ker(d\colon
L^i\to L^{i+1})$, $B^i(L)=image(d\colon
L^{i-1}\to L^{i})$, $H^i(L)=Z^i(L)/B^i(L)$.\\

\begin{definition}\label{1.3}
The Maurer-Cartan equation (also called the deformation
equation) of a DGLA $L$ is
\[ da+\frac{1}{2}[a,a]=0,\qquad a\in L^1.\]
\end{definition}

There is also an obvious notion of morphisms of DGLA's; every morphism
of DGLA induces a morphism between cohomology groups.  It is also
evident that morphisms of DGLA preserves solutions of the
Maurer-Cartan equation.\\
A {\it quasiisomorphism} is a morphism inducing isomorphisms in
cohomology.  Two DGLA's are {\it quasiisomorphic} if they are
equivalent under the equivalence relation generated by
quasiisomorphisms.\\

The cohomology of a DGLA is itself a differential graded Lie algebra
with the induced bracket and zero differential:\\

\begin{definition}\label{1.4}
A DGLA $L$ is called {\it Formal} if it is quasiisomorphic to its
cohomology DGLA $H^*(L)$.
\end{definition}

{\sc Exercise:} Let $D\colon L\to L$ be a derivation, then the kernel
of $D$ is a graded Lie subalgebra.\\

{\sc Exercise:} Let $L=\oplus L^i$ be a DGLA and $a\in L^i$.  Prove
that:\begin{enumerate}

\item If $i$ is even then $[a,a]=0$.

\item If $i$ is odd then $[a,[a,b]]=\dfrac{1}{2}[[a,a],b]$ for every
$b\in L$ and $[[a,a],a]=0$.\end{enumerate}

In nature there exists several examples of DGLA, most of which govern
deformation problems via Maurer-Cartan equation; here we present some of
the most interesting ones.\\

\begin{example}\label{1.5}
Every Lie algebra is a DGLA concentrated in degree 0.\end{example}

\begin{example}\label{1.6} Let $A=\oplus A^i$ be an associative  graded-commutative
$\K$-algebra (this means that $ab=(-1)^{\overline{a}\overline{b}}ba$ for $a,b$
homogeneous) and $L=\oplus L^i$ a DGLA. Then $L\otimes_{\K}A$ has a
natural structure of DGLA by setting:\\
\[ (L\otimes_{\K}A)^n=\oplus_i (L^{i}\otimes_{\K}A^{n-i}),\quad
d(x\otimes a)=dx\otimes a,\quad
[x\otimes a,y\otimes b]=(-1)^{\overline{a}\overline{y}}[x,y]\otimes ab.\]
\end{example}

In the next examples, when
$M$ is a complex variety and $E$ a holomorphic vector bundle on $M$,
we denote by $\sA^{p,q}(E)$ the sheaf of $C^{\infty}$ $(p,q)$-differential
forms with values in $E$.

\begin{example}\label{1.7} Let $E$ be a holomorphic vector bundle on a
complex variety $M$.  We have a DGLA $L=\oplus \Gamma(M,
\sA^{0,p}(End(E)))[-p]$ with the Dolbeault differential and the
natural bracket.  More precisely if $e,g$ are local holomorphic
sections of $End(E)$ and $\phi,\psi$ differential forms we define
$d(\phi e)=(\overline{\partial}\phi)e$, $[\phi e,\psi
g]=\phi\wedge\psi[e,g]$.\end{example}

\begin{theorem}\label{1.8} If $M$ is K\"{a}hler and $E$ is a flat unitary stable
bundle (e.g. $E=\Oh_M$) then the DGLA of \ref{1.7} is
formal.\end{theorem}

\begin{proof} \cite[p. 92]{GM1}. \end{proof}

\begin{example}\label{1.9} Let
$T_M$ be the holomorphic tangent bundle of a complex variety $M$. The
Kodaira-Spencer DGLA is defined as $KS(M)=\oplus \Gamma(M,\sA^{0,p}(T_M))[-p]$
with the
Dolbeault differential; if $z_1,\ldots,z_n$ are local holomorphic
coordinates we have $[\phi d\overline{z}_I, \psi
d\overline{z}_J]=[\phi,\psi]d\overline{z}_I\wedge d\overline{z}_J$ for
$\phi,\psi\in \sA^{0,0}(T_M)$, $I,J\in
\bigwedge^*\{1,\ldots,n\}$.\end{example}

\begin{example}\label{1.10} 
Let $\K=\R$ and $X$ be a smooth differentiable
variety. The algebra of polyvector fields is given by
\[ T_{Poly}(X)=\oplus \Gamma(X,\bigwedge^{n+1}T_X)[-n],\qquad n\ge 
-1\]
with zero differential and the \emph{Schouten-Nijenhuis} bracket defined
in the following way (cf. \ref{1.12}):\\
For every open subset $U\subset X$, every function
$h\in T^{-1}_{Poly}(U)=C^{\infty}(U)$ and every vector fields
$\xi_0,\ldots,\xi_n,\zeta_0,\ldots,\zeta_m\in \Gamma(U, T_U)$ we set
\[ [\xi_0\wedge\cdots\wedge\xi_n,h]=\sum_{i=0}^n(-1)^{n-i}
\xi_i(h)\xi_0\wedge\cdots\wedge
\widehat{\xi_i}\wedge\cdots\wedge\xi_n\]
\begin{multline*} 
[\xi_0\wedge\cdots\wedge\xi_n,\zeta_0\wedge\cdots\wedge\zeta_m]=\\
=\sum_{i=0}^n\sum_{j=0}^m
(-1)^{i+j}[\xi_i,\zeta_j]\wedge\xi_0\wedge\cdots\wedge
\widehat{\xi_i}\wedge\cdots\wedge\xi_n\wedge\zeta_0\wedge\cdots\wedge
\widehat{\zeta_j}\wedge\cdots\wedge\zeta_m.\end{multline*}
\end{example}

\begin{example}\label{1.11}
Let $A$ be an associative $\K$-algebra, the DGLA of Hochschild
cochains is defined by
\[ G=\oplus \Hom_{\K}(A^{\otimes (n+1)},A)[-n]\]
where by $\Hom_{\K}$ we mean homomorphisms of $\K$-vector spaces. The
differential is the usual differential of Hochschild cohomology:
$\phi\in G^{n}$
\[  (d\phi)(a_0\otimes\cdots\otimes a_{n+1})
=a_0\phi(a_0\otimes\cdots\otimes a_{n+1})+\]
\[
+(-1)^n\phi(a_0\otimes\cdots\otimes a_n)a_{n+1}+\sum_{i=0}^n(-1)^i
\phi(a_0\otimes\cdots\otimes a_i a_{i+1}\otimes\cdots\otimes  
a_{n+1}).\]
The bracket is the \emph{Gerstenhaber} one:
\[ [\phi,\psi]=\phi\circ\psi-(-1)^{\overline{\phi}\overline{\psi}}\psi\circ 
\phi\]
where the (non-associative) product $\circ$ is defined, for
$\phi\in G^n$ and $\psi\in G^m$, by the formula
\begin{multline*} 
\phi\circ \psi(a_0\otimes\cdots\otimes a_{n+m})=\\
=\sum_{i=0}^{n}(-1)^{im}\phi(a_0\otimes\cdots\otimes a_{i-1}\otimes
\psi(a_i\otimes\cdots\otimes a_{i+m})\otimes a_{i+m+1}\otimes\cdots\otimes 
a_{n+m}).\end{multline*}
\end{example}

\begin{example}\label{1.12} Let $A$ be a commutative $\K$-algebra and $V\subset
Der_{\K}(A,A)$ an $A$-submodule such that $[V,V]\subset V$.\\
We can define a DGLA $L=\oplus_{n\ge -1} L^n$, where
\[ L^{-1}=A,\quad L^0=V,\ldots,  L^n=\bigwedge_A^{n+1}V,\ldots\] 
with zero differential and the bracket uniquely characterized by the
properties:\begin{enumerate}
\item $[,]\colon L^0\times L^0\to L^0$ is the usual bracket on $V$.

\item If $a\in L^0$, $f\in L^{-1}=A$ then $[a,f]=a(f)\in A$.

\item For every $a\in L^n$, $b\in L^m$, $c\in L^h$
\[ [a,b\wedge c]=[a,b]\wedge c+(-1)^{n(m+1)}b\wedge [a,c],\quad
[a\wedge b,c]=a\wedge [b,c]+(-1)^{h(m+1)}[a,c]\wedge b.\] 
\end{enumerate}
\end{example}

{\sc Exercise:} Prove that the bracket in \ref{1.12} is well
defined. (the unicity is obvious, the existence is easy when $V$
is free; in the general case it is convenient to think $L^n$ as
the quotient of the $\K$-vector space generated by
$\xi_0\wedge\cdots\wedge\xi_n$, with $\xi_i\in
V$, by the subspace generated by the skewsymmetric $A$-multilinear relation).\\
Prove moreover that for every
$h\in L^{-1}$, $\xi_0,\ldots,\xi_n,\zeta_0,\ldots,\zeta_m\in L^0$ we have (cf. \ref{1.10})
\[ [\xi_0\wedge\cdots\wedge\xi_n,h]=
\sum_{i=0}^n(-1)^{n-i}\xi_i(h)\xi_0\wedge\cdots\wedge
\widehat{\xi_i}\wedge\cdots\wedge\xi_n\]
\begin{multline*} 
[\xi_0\wedge\cdots\wedge\xi_n,\zeta_0\wedge\cdots\wedge\zeta_m]=\\
=\sum_{i=0}^n\sum_{j=0}^m
(-1)^{i+j}[\xi_i,\zeta_j]\wedge\xi_0\wedge\cdots\wedge
\widehat{\xi_i}\wedge\cdots\wedge\xi_n\wedge\zeta_0\wedge\cdots\wedge
\widehat{\zeta_j}\wedge\cdots\wedge\zeta_m.\end{multline*}

For a better understanding of some of next topics it is useful to
consider the following functorial construction. Given
a DGLA $(L, [,],d)$ we can construct a new DGLA $(L_d, [,]_d, d_d)$ by
setting
$L^i_d=L^i$ for every $i\not=1$, $L^1_d=L^1\oplus\K d$ with the bracket
defined by extending by bilinearity the relation
$[a+vd,b+wd]_d=[a,b]+vd(b)+(-1)^{\overline{a}}wd(a)$ and differential
$d_d(a+vd)=[d,a+vd]_d=d(a)$.\\
The natural inclusion $L\subset L_d$ is a morphism of DGLA; in the manage
of Maurer-Cartan equation it is convenient to consider the affine
embedding $\phi\colon L\to L_d$, $\phi(a)=a+d$.\\
For $a\in L^1$ we have
\[ d(a)+\frac{1}{2}[a,a]=0\quad\iff\quad [\phi(a),\phi(a)]_d=0.\]

\begin{example}\label{1.13} (Small variations of almost complex
structures).\\
Let $M$ be a compact complex variety, $A^{p,q}$
the vector bundle of  differential
forms of type $(p,q)$,
$\sA^{p,q}$ the fine sheaf of its smooth sections
and $KS(M)=\oplus \Gamma(M,\sA^{0,p}(T_M))[-p]$ the
Kodaira-Spencer algebra of $M$.
Note that $\Gamma(M,\sA^{0,p}(T_M))$ is exactly the set of morphisms of vector
bundles $A^{1,0}\to A^{0,p}$.
The choice of a
hermitian metric on $M$ induce a structure of pre-Hilbert space to
$\Gamma(M,\sA^{0,p}(T_M))$.
Given a sufficiently small section $\alpha\in \Gamma(M, \sA^{0,1}(T_M))$ we can
consider the perturbed differential
$\overline{\partial}_\alpha=\overline{\partial}+\alpha\partial\colon
\sA^{0,p}\to \sA^{0,p+1}$; more precisely for $U\subset M$ open subset
with holomorphic coordinates $z_1,\ldots,z_n$ and $f\in C^{\infty}(U)$
we have
$\overline{\partial}_\alpha(fd\overline{z}_I)=(\overline{\partial}f
+\alpha(\partial f))\wedge d\overline{z}_I$.\\
To $\alpha$ we can associate also a direct sum decomposition
\[ A^{1,0}_\alpha\oplus A^{0,1}_\alpha=A^{1,0}\oplus A^{0,1}\]
where $A^{1,0}_{\alpha}$ is the graph of $-\alpha$ and
$A^{0,1}_\alpha=\overline{A^{0,1}_\alpha}$. Therefore we can consider $\alpha$
as a small variation of the almost complex structure, as easy computation
(exercise) show that:\begin{enumerate}

\item The sheaf of almost holomorphic functions of the almost
complex structure $A^{i,j}_\alpha$ is exactly the kernel of
$\overline{\partial}_\alpha\colon \sA^{0,0}\to \sA^{0,1}$.

\item $\alpha$ satisfies the
Maurer-Cartan equation in the Kodaira-Spencer algebra
if and only if $\overline{\partial}_\alpha^2=0$.

\item By Newlander-Nirenberg theorem the almost complex structure
$A^{i,j}_\alpha$ is integrable if and only if $\alpha$ satisfies the
Maurer-Cartan equation in $KS(M)$.
\end{enumerate}
\end{example}

Let's now introduce the notion of gauge action on a DGLA:\\
There exists a functor $exp$ from the category
of nilpotent Lie algebras (i.e. with descending central series definitively =0)
to the category of groups. For every
nilpotent Lie algebra $N$ there exists a natural bijection $e\colon
N\to exp(N)$ satisfying the following properties:\begin{enumerate}

\item $e^ae^b=e^{(a+b)}$ if $[a,b]=0$.
More generally $e^ae^b=e^{(a*b)}$ where $a*b\in N$ is given by
the Campbell-Baker-Hausdorff formula (cf \cite{ea}, \cite{Ho}, \cite{Ja}).

\item For every vector space $V$ and every
homomorphism of Lie algebras $\rho\colon N\to End(V)$ such that $\rho(N)$
is a nilpotent subalgebra, the morphism
\[exp(\rho)\colon exp(N)\to Aut(V),\qquad exp(\rho)(e^a)=e^{\rho(a)},\]
is a
homomorphism of groups (here $e^{\rho(a)}$ denotes the usual exponential of
endomorphisms).

\item If $\rho\colon N\to End(V)=P$ is a representation of $N$ as above and
$ad_\rho\colon N\to End(P)$ is the adjoint representation, then
for every $a\in N$, $f\in End(V)$
\[ exp(ad_\rho)(e^a)f=e^{\rho(a)}\circ f\circ e^{-\rho(a)}.\]
\end{enumerate}

{\sc Exercise:} Prove the above items.\\

\begin{lem}\label{1.14} Let $V,W$ be vector spaces and $\rho\colon N\to
End(V)$, $\eta\colon N\to End(W)$ representation of a nilpotent Lie
algebra $N$.\\
Let $b\colon V\times V\to W$ be a bilinear symmetric form
such that $2b(v,\rho(n)v)=\eta(n)b(v,v)$ for every $n\in N$, $v\in V$
and
$q(v)=b(v,v)$ the associated quadratic form.\\
Then the cone $Z=q^{-1}(0)$ is invariant under the exponential action
$exp(\rho)$.\end{lem}

\begin{proof} Let $n\in N$ be a fixed element, for every $v\in V$ define the
polynomial function $F_v\colon\K\to W$ by
\[ F_v(t)=exp(\eta(-tn))q(exp(\rho(tn))v).\]
For every $s,t\in\K$, if $u=exp(\rho(sn))v$ then
\[ F_v(t+s)=e^{\eta(-sn)}F_w(t),\quad
\frac{\partial F_v}{\partial t}(0)=-\eta(n)q(v)+2b(v,\rho(n)v)=0\]
\[ \frac{\partial F_v}{\partial t}(s)=
e^{\eta(-sn)}\frac{\partial F_w}{\partial t}(0)=0\]
As the field $\K$ has characteristic 0 every function $F_v$ is
constant, in particular for every $exp(\eta)$ invariant subspace
$T\subset W$ the space $q^{-1}(T)$ is $exp(\rho)$ invariant.\end{proof}

\begin{cor}\label{1.15}
 Let $L$ be a nilpotent DGLA, then the quadratic cone
\[Z=\{a\in L^1\mid [v,v]=0\}\]
is stable under the exponential of the
adjoint action of $L^0$.\end{cor}

\begin{proof}
By Jacobi identity $[v,[a,v]]=-[v,[v,a]]=\dfrac{1}{2}[a,[v,v]]=0$
and we can apply \ref{1.14}.\end{proof}

{\sc Exercise:} In the notation of \ref{1.14},
if $\K$ is algebraically closed and $q^{-1}(0)$ is not a double plane
then \ref{1.14} holds under the weaker assumption
$b(v,\rho(n)v)=0$ for every $n\in N$, $v\in Z$.\\

Again if $L$ is a DGLA with $L^0$ nilpotent, the adjoint action of $L^0$
over $L^1_d$ preserves the cone $[v,v]_d=0$ and the affine hyperplane
$\{v+d\mid v\in L^1\}$; this gives an action of $exp(L^0)$ over $L^1$
preserving the solution of the Maurer-Cartan equation. The infinitesimal
generator of this action is given by
$L^0\ni a\mapsto [a, v+d]_d=[a,v]-d(a)\in T_vL^1$.\\

\begin{remark} It is often convenient to think the elements of $L_d$ as
operators
on a $\Z$-graded vector space $V=\oplus V^i$, this means that
the map $\rho\colon L_d\to End(V)=P$ is a morphism of DGLA.\\
For example the elements of the extended Kodaira-Spencer algebra
\[ KS(M)_{\partial}=\oplus \Gamma(M,\sA^{0,p}(T_M))[-p]\oplus\C\overline{\partial}\]
act in a natural way on the graded vector space
$V=\oplus \Gamma(M,\sA^{0,p})[-q]$.\\
If we call $ad_{\rho}\colon L_d\to End(P)$ and
$Ad_\rho\colon GL(L_d)\to GL(P)$ the
the adjoint actions
\[ ad_\rho(a)B=[\rho(a),B]=\rho(a)B-B\rho(a),\qquad Ad_\rho(e^a)B=
e^{\rho(a)}B e^{-\rho(a)}\]
by the properties of $exp$ we have $exp(ad_\rho(a))=Ad_\rho(exp(a))$.\\
It is natural to consider two operator $A,B\in End(V)$ gauge equivalent if there
exists $a\in L^0$ such that $A=e^{\rho(a)} B  e^{-\rho(a)}$.\\
Therefore if $x,y\in L^1_d$ are gauge equivalent then $\rho(x)$ is gauge
equivalent to $\rho(y)$, independently from the particular representation.
\end{remark}

\begin{example}\label{1.16} A pair $(L,d)$
is a \emph{Polarized graded algebra} (PGA) if $L=\oplus L^i$, $i\in \Z$, is
a graded associative $\K$-algebra, $d\in L^1$ and $d^2=0$.\\
The typical example is, for a fixed
differential graded vector space $(V,d)$, $d\colon V\to V[1]$,
the algebra $End(V)=\oplus L^i$, where $L^i$
is the space of morphisms $V\to V[i]$ in the category
$Graded^{\K}$.\\
Every PGA admits  natural structures of DGA (differential graded algebra)
with differential $\delta(a)=da-(-1)^{\overline{a}}ad$ and  DGLA with bracket
$[a,b]=ab-(-1)^{\overline{a}\overline{b}}ba$.\\
If $L^0$ is a nilpotent $\K$-algebra then we can define the gauge
action over $L^1$ in a characteristic free way. First we define a group
structure on $L^0$ by setting $g*h=g+h+gh$ and the gauge action
\[ g*v=v+\sum_{i=0}^\infty (-1)^i[g,d+v]g^i=v+\sum_{i=0}^\infty
(-1)^i([g,v]-\delta(g))g^i.\]
\end{example}

{\sc Exercise:} Prove that in characteristic 0 the action of \ref{1.16} is
equivalent to the usual gauge action.\\

{\sc Exercise:} Let $L$ be a DGLA with differential $d$,
then the universal enveloping algebra
of $L_d$ is a polarized graded algebra.\\

In some cases (as in  \cite[p. 9]{K1})
it is useful to describe a DGLA as a suitable
subset of a concrete PGA.\\

\bigskip

\section{Quickstart guide to functors of Artin rings}
\label{sec.2}

In this section $\K$ is  a fixed field of arbitrary characteristic.
We denote by:\\
$\widehat{Art}$ the category of local complete Noetherian
rings with residue field $\K$.\\
$Art\subset \widehat{Art}$ the full subcategory of Artinian rings.\\
For a given $S\in \widehat{Art}$,  $Art_S$ is the category of Artinian
$S$-algebras with residue field $\K$.\\
$\widehat{Art}_S$ the category of local complete Noetherian
$S$-algebras with residue field $\K$.\\
$Set$ the category of sets (in a fixed universe).
$*$ is a fixed 1-point set.\\
$Grp$ the category of groups. ($\emptyset$ is not a group)\\

\begin{definition}\label{2.1}
A \emph{functor of Artin rings} is a covariant functor
$F\colon Art_S\to Set$, $S\in \widehat{Art}$ such that $F(\K)=*$.
\end{definition}

The main interest to functors of Artin rings comes from deformation
theory and moduli problems; from this point of view the notion of
prorepresentability is one of the most important.\\

Given $R\in \widehat{Art}_S$ we define a functor $h_R\colon Art_S\to Set$ by
setting $h_R(A)$ as the set of $S$-algebra homomorphisms $R\to A$.

\begin{definition}\label{2.2}
A functor $F\colon Art_S\to Set$ is {\it prorepresentable}
if it is isomorphic to $h_R$ for some $R\in
\widehat{Art}_S$.
\end{definition}

The functors of Artin rings $F\colon Art_S\to Set$, with their natural
transformation form a category denoted by $Fun_S$.\\
We left as an exercise
to prove that the Yoneda functor $\widehat{Art}_S^0\to Fun_S$, $R\to h_R$, is
fully faithful.\\

A necessary  condition for a functor $F$ to be
prorepresentable is homogeneity. We first note that on $Art_S$ there
exist fibred products
\begin{equation}\label{2.3}
\begin{array}{ccc}
B\times_AC&\mapor{}&C\\
\mapver{}&&\mapver{}\\
B&\mapor{}&A\end{array}
\end{equation}

Applying a functor  $F\in Fun_S$ to the cartesian diagram (\ref{2.3}) we get a
map
\[ \eta \colon F(B\times_AC)\to F(B)\times_{F(A)}F(C).\]

\begin{definition}\label{2.4}
The functor $F$ is {\it homogeneous} if $\eta$ is an
isomorphism whenever $B\to A$ is surjective.\end{definition}

Since the diagram $(\ref{2.3})$ is cartesian,
every prorepresentable functor is homogeneous.

\begin{example}\label{2.5}
Other examples of homogeneous functors are:\\
\ref{2.5}.0) The trivial functor $F(A)=*$ For every $A\in Art_S$.\\
\ref{2.5}.1) Let $M$ be a flat $S$-module, define $\widehat{M}\colon Art_S\to Set$ as
$\widehat{M}(A)=M\otimes_S m_A$, where $m_A\subset A$ is the maximal ideal.\\
\ref{2.5}.2) Assume $char\K=0$ and let $L^0$ be a Lie Algebra over $\K$.
We can define a group functor $exp(L^0)\colon Art_{\K}\to Grp$ by
setting $exp(L^0)(A)=exp(L^0\otimes m_A)$ the exponential of the
nilpotent Lie algebra $L^0\otimes m_A$.\\
\ref{2.5}.3) Let $X\to Spec(S)$ be a flat scheme over $S$, we can define
$Aut(X/S)\colon Art_S\to Grp$ by setting $Aut(X/S)(A)$ as the group of
automorphisms of $X_A=X\times_{Spec(S)}Spec(A)$, commuting with the projection
$X_A\to Spec(A)$ which are the identity on $X_{\K}$.
\end{example}

{\sc Exercise:} Prove that the above functors $\widehat{M}$, $exp(L^0)$ and
$Aut(X/S)$ are homogeneous. Show also by an example that if $M$ is not
$S$-flat then the functor $\widehat{M}$ is not homogeneous in general.\\

\begin{definition}\label{2.6}
A functor $F$ is called a \emph{deformation functor}
if:\begin{enumerate}
\item  $\eta$ is surjective whenever $B\to A$ is surjective.

\item $\eta$ is an isomorphism whenever $A=\K$.
\end{enumerate}
\end{definition}

The name comes from the fact that most  functors arising in
deformation theory  are deformation functors.\\
{\sc Exercise:}
Let $X$ be a scheme over $\K$, the deformation functor of $X$ is
defined as \\
$Def_X\colon Art_\K\to Set$, where $Def_X(A)$ is the set of
isomorphism classes of commutative cartesian diagrams
\[ \begin{array}{ccc}
X&\mapor{i}&X_A\\
\mapver{p}&&\mapver{p_A}\\
Spec(\K)&\mapor{}&Spec(A).\end{array}\]
with $i$ closed embedding and $p_A$ flat morphism.
Prove that $Def_X$ is a deformation functor.\medskip\\

We denote  $\K[\epsilon]=\K\oplus \K\epsilon=\K[t]/(t^2)$ with the trivial
structure of $S$-algebra given by $S\to\K\to\K[\epsilon]$. More generally
by $\epsilon$ and $\epsilon_i$ we will always mean indeterminates annihilated by the
maximal ideal, and in particular of square zero (e.g., the algebra $\K[\epsilon]$
has dimension $2$ and $\K[\epsilon_1,\epsilon_2]$ has dimension $3$ as a $\K$-vector
space).

\begin{prop}[Schlessinger, \cite{Sch}]\label{2.7}
    Let $F$ be a deformation functor, the set
$t_F=F(\K[\epsilon])$ has a natural structure of $\K$-vector space.
If $\phi\colon F\to G$ is a morphism of deformation functors the map
$\phi\colon t_F\to t_G$ is linear.\end{prop}

\begin{proof}  Let $\alpha\in\K$, the scalar multiplication by $\alpha$ is
induced by the morphism in $Art_S$,
$\widehat{\alpha}\colon\K[\epsilon]\to\K[\epsilon]$ given by
$\widehat{\alpha}(a+b\epsilon)=a+\alpha b\epsilon$. Since
$t_F\times t_F=F(\K[\epsilon_1]\times_{\K}\K[\epsilon_2])$ the sum is induced by
the map $\K[\epsilon_1]\times_{\K}\K[\epsilon_2]=
\K[\epsilon_1,\epsilon_2]\to\K[\epsilon]$
defined by $a+b\epsilon_1+c\epsilon_2\to a+(b+c)\epsilon$.\end{proof}

{\sc Exercise:} Prove that if $R\in \widehat{Art}_{\K}$,
$F=h_R$ then $t_F$ is isomorphic to the
Zariski tangent space of $Spec(R)$;
more generally if $R\in \widehat{Art}_{S}$ then the tangent space of
$h_R\in Fun_S$ is equal to the $\K$-dual of $m_R/(m^2_R+m_SR)$.\\

{\sc Exercise:} Let $A\in Art_S$, $F\in Fun_S$. Prove that there exists a
natural bijection between $F(A)$ and morphisms $Mor(h_A,F)$.

\begin{definition}\label{2.8}
A morphism $\phi\colon F\to G$ in the category $Fun_S$ is
called:\begin{enumerate}
\item \emph{unramified} if $\phi\colon t_F\to t_G$ is injective.

\item \emph{smooth} if for every surjection  $B\to A$ in
$Art_S$ the map
$F(B)\to G(B)\times_{G(A)}F(A)$ is also
surjective.

\item \emph{\'etale} if it is smooth and unramified.
\end{enumerate}
\end{definition}

{\sc Exercise:} If $\phi\colon F\to G$ is smooth then $F(A)\to G(A)$ is
surjective for every $A$.\\

{\sc Exercise:} Let $\phi\colon F\to G$, $\psi\colon H\to G$ be
morphisms of deformation functors. If $\phi$ is smooth and $H$ is
prorepresentable then there exists a morphism $\tau\colon H\to F$ such
that $\psi=\phi\tau$.\\

\begin{definition}\label{2.9}
A functor $F$ is {\it smooth} if the morphism $F\to *$ is
smooth, i.e. if $F(A)\to F(B)$ is surjective for every surjective
morphism of $S$-algebras $A\to B$.\end{definition}

\begin{lem}\label{2.10}
A prorepresentable functor $h_R$, $R\in \widehat{Art}_S$, is
smooth if and only if $R=S[[x_1,...,x_n]]$.\end{lem}

\begin{proof} Exercise, cf. \cite{Sch}.\end{proof}

{\sc Exercise:} Let $X$ be a scheme over $\K$. Prove that
if for every $A\in Art_\K$ and every deformation
$X_A\to Spec(A)$ the functor $Aut(X/A)$ is smooth then $Def_X$ is
homogeneous.\medskip\\

\begin{lem}\label{2.11}
Let $\phi\colon F\to G$ be an unramified morphism of
deformation functors. If $G$ is homogeneous then $\phi\colon F(A)\to G(A)$ is
injective for every $A$.\end{lem}

\begin{proof}  We prove by induction on the length
of $A$ that $F(A)\to G(A)$ is injective.\\
Let $\epsilon\in A$ be an element such that $\epsilon m_A=0$ and $B=A/(\epsilon)$; by
induction $F(B)\to G(B)$ is injective.
We note that
\[ A\times_{\K}\K[\epsilon]\to A\times_B A,\qquad
(a, \overline{a}+\alpha\epsilon)\mapsto (a, a+\alpha\epsilon)\] 
is an isomorphism of $S$-algebras, in particular for every deformation
(resp. homogeneous)
functor $F$ there exists a natural surjective (resp. bijective) map
\[ v\colon F(A)\times t_F\to F(A)\times_{F(B)}F(A)\]
such  that $v(F(A)\times\{0\})=\Delta$ is the diagonal.
Let $\xi,\eta\in F(A)$ such that $\phi(\xi)=\phi(\eta)$; since $\phi$ is
injective over $F(B)$ the pair $(\xi,\eta)$ belongs to the fibred product
$F(A)\times_{F(B)}F(A)$ and there exists $h\in t_F$ such that
$v(\xi,h)=(\xi,\eta)$. Therefore
$v(\phi(\xi),\phi(h))=(\phi(\xi),\phi(\eta))\in \Delta$ and since $G$ is
homogeneous $\phi(h)=0$. By the definition of unramified morphism
$\phi\colon t_F\to t_G$ is injective and then $h=0$, $\xi=\eta$.\end{proof}

\begin{cor}\label{2.12} Let $\phi\colon F\to G$ be an \'etale morphism of
deformation functors. If $G$ is homogeneous then $\phi$ is an
isomorphism.\end{cor}

\begin{proof} Evident.\end{proof}

Let $\sim$ the equivalence relation on the category  of deformation
functors generated by the \'etale morphisms. It is an easy consequence (left
as exercise) of the above results
that $h_R\sim h_T$ if and only if $R,T$
are isomorphic $S$-algebras.\\

By a small extension $e$ in $Art_S$ we mean an exact sequence
\[ e\colon\quad 0\mapor{}M\mapor{}B\mapor{\phi}A\mapor{}0\]
where $\phi$ is a morphism in $Art_S$ and $M$ is an ideal of $B$ annihilated
by the maximal ideal $m_B$. In particular $M$ is a finite dimensional
vector space over $B/m_B=\K$.\\

{\sc Exercise:} Let $0\mapor{}M\mapor{}B\mapor{\phi}A\mapor{}0$ be a
small extension and $F$ a deformation functor. Then there exists
a natural transitive action of $F(\K\oplus M)=t_F\otimes M$ on every
nonempty fibre of $F(A)\mapor{\phi}F(B)$. (Hint: look at the proof of
\ref{2.11}).\medskip\\

\begin{definition}\label{2.13}
Let $F$ be a functor of Artin rings; an {\it obstruction theory}
$(V,v_e)$ for
$F$ is the data of a $\K$-vector space $V$, called {\it obstruction space},
and for every small extension in $Art_S$
\[ e\colon\quad 0\mapor{}M\mapor{}B\mapor{}A\mapor{}0\] 
of an {\it obstruction map}
$v_e\colon F(A)\to V\otimes_{\K}M$ satisfying the following
properties:\begin{enumerate}
\item If  $\xi\in F(A)$ can be lifted to $F(B)$ then
$v_e(a)=0$.

\item (base change) For every morphism $\alpha\colon e_1\to e_2$
of small extension, i.e. for
every commutative diagram
\[ \begin{array}{cccccccccc}
e_1\colon\quad&0&\mapor{}&M_1&\mapor{}&B_1&\mapor{}&A_1&\mapor{}&0\\
&&&\mapver{\alpha_M}&&\mapver{\alpha_B}&&\mapver{\alpha_A}&&\\
e_2\colon\quad&0&\mapor{}&M_2&\mapor{}&B_2&\mapor{}&A_2&\mapor{}&0.\end{array}
\]
we have
$v_{e_2}(\alpha_A(a))=(Id_V\otimes\alpha_M)(v_{e_1}(a))$ for every $a\in
F(A_1)$.
\end{enumerate}
\end{definition}

{\sc Exercise:} If $F$ is smooth then all the obstruction maps are
trivial.

\begin{definition}\label{2.14}
An obstruction theory $(V,v_e)$ for $F$ is called {\it
complete} if the converse of item i) in \ref{2.13} holds; i.e. the lifting
exists if and only if the obstruction vanish.\end{definition}

Clearly if $F$ admits a complete obstruction theory then it admits infinitely
ones; it is in fact sufficient to embed $V$ in a bigger vector space. One
of the main interest is to look for the ``smallest'' complete obstruction
theory.\\

\begin{definition}\label{2.15}
A morphism of obstruction theories $(V,v_e)\to (W,w_e)$ is a
linear map $\theta\colon V\to W$ such that $w_e=\theta v_e$ for every
small extension $e$. An obstruction theory $(O_F,ob_e)$ for $F$
is called {\it universal} if for every obstruction theory $(V,v_e)$ there
exist an unique morphism $(O_F,ob_e)\to (V,v_e)$.\end{definition}

If a universal obstruction theory exists then it is unique up to
isomorphisms. An important result is

\begin{theorem}\label{2.16} Let $F$ be a deformation functor, then there exists
an universal obstruction theory $(O_F,ob_e)$ for $F$.
Moreover the universal obstruction theory is complete and
every element of the vector space $O_F$ is of the form
$ob_e(\xi)$ for some principal extension
\[ e\colon\quad 0\mapor{}\K\mapor{}B\mapor{}A\mapor{}0\]
and some $\xi\in F(A)$.\end{theorem}

\begin{proof} This is quite long and not easy. The interested reader can found
a proof in  \cite{FM1}.\end{proof}

We note that if $F$ is not a deformation functor then in general $F$ doesn't
have any complete obstruction theory even if $F$ satisfies
Schlessinger's conditions H1, H2, H3 of \cite{Sch}.\\

\begin{example}\label{2.17} (The primary obstruction map, $char\K\not=2$).\\
Let $(V,v_e)$ be a complete  obstruction theory for a deformation
functor $F$ and let $[,]\colon t_F\times t_F\to V$ be the obstruction map
associated to the small extension
\[ 0\mapor{}\K\mapor{xy}\K[x,y]/(x^2,y^2)
\mapor{}\K[x,y]/(x^2,xy,y^2)\mapor{}0\] 
(Note that $\K[x,y]/(x^2,xy,y^2)=\K[x]/(x^2)\times_{\K}\K[y]/(y^2)$).
Then $[,]$ is a symmetric bilinear map (Exercise).\\
The substitution $\alpha(t)=x+y$ gives a morphism of small extensions
\[ \begin{array}{ccccccccc}
0&\mapor{}&\K&\mapor{t^2}&\K[t]/(t^3)&\mapor{}&\K[t]/(t^2)
&\mapor{}&0\\
&&\mapver{2}&&\mapver{\alpha}&&\mapver{\alpha}&&\\
0&\mapor{}&\K&\mapor{xy}&\K[x,y]/(x^2,y^2)&
\mapor{}&\K[x,y]/(x^2,xy,y^2)&\mapor{}&0.\end{array}\] 
From this and base change axiom it
follows that the obstruction of lifting $\xi\in t_F$ to
$F(\K[t]/(t^3))$ is equal to $\dfrac{1}{2}[\xi,\xi]$ ``the quadratic
part of Maurer-Cartan equation ''!!\end{example}

Let $\phi\colon F\to G$ be a morphism of deformation functors and
$(V,v_e)$, $(W,w_e)$ obstruction theories for $F$ and $G$
respectively; a linear map $\phi'\colon V\to W$ is {\it compatible} with
$\phi$ if $w_e\phi=\phi' v_e$ for every small extension $e$.\\

\begin{prop}[Standard smoothness criterion]\label{2.18}
Let $\nu\colon F\to G$ be a morphism of deformation
functors and $(V,v_e)\mapor{\nu'}(W,w_e)$ a compatible morphism between
obstruction theories. If $(V,v_e)$ is complete,
$V\mapor{\phi'}W$ injective and $t_F\to t_G$ surjective then $\phi$ is
smooth.\end{prop}

\begin{proof}
Let $e\colon 0\mapor{}\K\mapor{}B\mapor{}A\mapor{}0$ be a small
extension and let $(a,b')\in F(A)\times_{G(A)}G(B)$;
let $a'\in G(A)$ be the common image of $a$ and $b'$. Then
$w_e(a')=0$, as $a'$ lifts to $G(B)$, hence $v_e(a)=0$ by injectivity of
$\nu'$.
Therefore $a$ lifts to some $b\in F(B)$. In general $b''=\nu(b)$ is not
equal to $b'$. However, $(b'',b')\in G(B)\times_{G(A)}G(B)$
and therefore $b''$
differs from $b'$ by the action of an element $v\in t_G$ ($v$ need not be
unique). As $t_F\to t_G$ is surjective, $v$ lifts to a $w\in t_F$; acting
with $w$ on $b$ produces a lifting of $a$ which maps to $b'$, as required.
\end{proof}

The composition $(V,v_e\phi)$ is an
obstruction theory for $F$ and therefore there exists a unique compatible
map $O_F\mapor{o(\phi)} V$.
In particular this happens for the universal obstruction
of $G$. Therefore every morphism of deformation functors $F\to G$
induces linear morphisms both in tangent $t_F\to t_G$ and
obstruction spaces $O_F\to O_G$. As an immediate consequence of \ref{2.18} we
have

\begin{prop}\label{2.19}
A morphism of deformation functors $\nu\colon F\to G$ is smooth
if and only if $t_F\to t_G$ is surjective and $o(\nu)\colon
O_F\to O_G$ is injective.
In particular $F$ is smooth if and only if $O_F=0$.\end{prop}

\begin{proof} One implication is contained in \ref{2.18}.
On the other side, if the morphism is smooth then $t_F\to t_G$ is
surjective; since every $x\in O_F$ is the obstruction to lifting some element
and $(O_G,ob_e)$ is complete we have $o(\nu)(x)=0$ if and only if $x=0$.
\end{proof}

{\sc Exercise:} Let $F\mapor{\phi}G\mapor{\psi}H$ be morphisms of
deformation functors. Prove that:\begin{enumerate}
\item If $\phi,\psi$ are smooth then the composition $\psi\phi$ is
smooth.

\item If $\psi\phi$, $\phi$ are smooth then $\psi$ is smooth.

\item If $\psi\phi$ is smooth and $t_F\to t_G$ is surjective then
$\phi$ is smooth.
\end{enumerate}

{\sc Exercise:} If $\phi\colon F\to G$ is smooth then $o(\phi)\colon
O_F\to O_G$ is an isomorphism.

\medskip
In most concrete cases it is very difficult to calculate the universal
obstruction space, while it is easy to describe complete obstruction
theories and compatible morphism between obstruction spaces.\\

Consider now the following situation.
$F\colon Art_{S}\to Set$ a deformation functor,
$G\colon Art_{S}\to Grp$ a group functor of Artin rings which is a
smooth deformation functor. (A theorem in \cite{FM1} asserts
that smoothness is automatic if $S$ is a field of characteristic 0).
\\
We assume that $G$ acts on $F$; this means that for every $A\in
Art_{S}$ there exists an action $G(A)\times F(A)\mapor{*} F(A)$, all these
actions must be compatible with morphisms  in $Art_{S}$.\\
In particular there exists an action $t_G\times t_F\mapor{*} t_F$, denote by
$\nu\colon t_G\to t_F$ the map $h\mapsto h*0$.\\

\begin{lem}\label{2.20} The map $\nu$ is linear and $t_G$ acts on $t_F$ by
translations, $h*v=\nu(h)+v$.\end{lem}

\begin{proof} Since the vector space structure on $t_F$ and $t_G$
is defined functorially by using morphisms in $Art_S$, it is easy to see
that for every $a,b\in t_F$, $g,h\in t_G$, $t\in\K$ we have
$(g+h)*(a+b)=(g*a)+(h*b)$, $t(g*a)=(tg)*(ta)$. Setting $a=b=0$ we get the
linearity of $\nu$ and setting $a=0$, $h=0$ we have
$g*b=(g*0)+(0*b)=\nu(g)+b$.\end{proof}

\begin{lem}\label{2.21} In the notation above the quotient functor
$D=F/G$ is a deformation functor,
$t_D=\Coker\nu$,  the projection $F\to D$ is smooth and for every
obstruction theory $(V,v_e)$ of $F$ the group functor
$G$ acts trivially on the obstruction maps $v_e$. In particular the
natural map $O_F\to O_D$ is an isomorphism.\end{lem}

\begin{proof} We left as exercise the (very easy) proof that $D$ is a
deformation functor and that $F\to D$ is a smooth morphism.
The statement about obstruction follow easily from \ref{2.16} and \ref{2.19}.
Since this result will be fundamental in  Section \ref{sec.3},
in order to make this lecture selfcontained we give here an
alternative proof which do not use the existence of $O_F,O_D$.\\
Let $a\in F(A)$, $g\in G(A)$ and
$e:0\mapor{}\K\mapor{j}B\mapor{}A\mapor{}0$ be a small extension. We
need to prove that $v_e(a)=v_e(g*a)$.\\
Let $\pi_i\colon A\times_{\K}A\to A$, $i=1,2$ be the projections,
$\delta\colon A\to
A\times_{\K}A$ be the diagonal and
\[ \nabla e\colon\quad 0\mapor{}\K\mapor{(j,0)}C\mapor{}A\times_{\K} 
A\mapor{}0\]
where $C$ is the quotient of $B\times_{\K}B$ by the ideal generated by
$(j,j)$.\\
Let $c=\delta(a)\in F(A\times_{\K}A)$ and let $\tau\in G(A\times_{\K}A)$
such that $\pi_1(\tau)=1$, $\pi_2(\tau)=g$, as $\delta$ lifts to a morphism
$A\to C$ and $G$ is smooth we have that $c$ and $\tau$ lift to $F(C)$ and
$G(C)$ respectively and therefore also $\tau * c$ lifts to $F(C)$.
Since $\pi_1({\tau}* c)={\pi_1(\tau)}*\pi_1(c)=a$,
$\pi_2({\tau}*c)={\pi_2(\tau)}*\pi_2(c)=g*a$.\\
On the other hand it is an easy consequence of the base change axiom
that $0=v_{\nabla e}(\tau *c)=v_e(\pi_1({\tau}* c))-v_e(\pi_2({\tau}*
c))=v_e(a)-v_e(g*a)$.\end{proof}

Let now $R\in Art_S$ e $a,b\in F(R)$, we define a functor
$Iso(a,b)\colon Art_R\to Set$ by setting
\[ Iso(a,b)(R\mapor{f}A)=\{g\in G(A)\mid g*f(a)=f(b)\}.\]

\begin{prop}\label{2.22} If $F$ is homogeneous then $Iso(a,b)$ is a deformation
functor with tangent space $\ker\nu$ and complete obstruction space
$\Coker\nu=t_D$.\end{prop}

\begin{proof}
For simplicity of notation let's denote $H=Iso(a,b)$. Assume it is
given a commutative diagram in $Art_{\K}$
\[ \begin{array}{ccc}
S&\mapor{\alpha}&A\\
\mapver{\beta}&&\mapver{\delta}\\
B&\mapor{\gamma}&C.\end{array}\]
with $\gamma$ surjective and let $(g_1,g_2)\in H(B)\times_{H(C)}H(A)$;
since $G$ is a deformation functor there exists $g\in G(B\times_C A)$
which lifts $(g_1,g_2)$. Let $a',b'$ be the images of $a,b$ on
$F(B\times_C A)$, the elements $g*a'$ and $b'$ have the same images in
$F(A), F((B)$, since $F$ is homogeneous it follows that $g\in H(B\times_C
A)$; this proves that $H$ is a deformation functor.\\
It is a trivial consequence of Lemma \ref{2.20} that $t_H=\ker\nu$.\\
Let consider a small extension in $Art_S$
\[ 0\mapor{}\K\mapor{\epsilon}A\mapor{}B\mapor{}0\]
and $g\in H(B)$, we want to determine the obstruction to lifting $g$ to
$H(A)$, by assumption $G$ is smooth and then there exists $g'\in G(A)$
lifting $g$. Let $a',b'\in F(A)$ be the images of $a,b$, then
\[ (g'*a',b')\in F(A)\times_{F(B)}F(A)=F(A\times_B A)=F(A)\times t_F.\]
It is easy to see that, changing $g'$ inside the liftings of $g$, the
projection over $t_F$ change by an element of the image of $\nu$,
therefore the natural projection of the pair $(g'*a',b')$ into
$\Coker\nu$ gives a complete obstruction (some details are left to
the reader).\end{proof}

The above construction of the obstruction theory for $Iso(a,b)$ is
natural, in particular if $G'$ is a smooth group
functors acting on a homogeneous functors $F'$ and $G\to G'$,
$\phi\colon F\to
F'$ are compatible morphisms then the induced morphism $\Coker\nu\to
\Coker\nu'$ is a morphism of obstruction spaces compatible with the
morphisms $Iso(a,b)\to Iso(\phi(a),\phi(b)$. We have therefore the following
corollary.

\begin{cor}\label{2.23}
In the above situation if $\ker\nu\to\ker\nu'$ is surjective
and $\Coker\nu\to \Coker\nu'$ is injective then $F/G\to F'/G'$ is
injective.\end{cor}
\begin{proof} (cf. \cite{Sch})
We need to prove that, given $a,b\in F(A)$, if
$Iso(\phi(a),\phi(b))(A)\not=\emptyset$ then $Iso(a,b)(A)\not=\emptyset$, by standard
criterion of smoothness the morphism $Iso(a,b)\to Iso(\phi(a),\phi(b))$
is smooth, in particular $Iso(a,b)(A)\to Iso(\phi(a),\phi(b))(A)$ is
surjective.\end{proof}

\begin{remark} The hypothesis $F$ homogeneous in Proposition
\ref{2.22} can be weakened by taking $F$ deformation functor and
assuming also the existence of
a complete relative obstruction theory for the morphism $G\to F$, $g\mapsto
g*F(\K)$, see \cite{FM1}. This last condition seems usually verified in concrete cases,
we only know counterexamples in positive characteristic.
\end{remark}

{\sc Exercise:} In the same notation of Corollary \ref{2.23}, if in addition
$G'$ is homogeneous
and $\ker\nu\to\ker\nu'$ is bijective then $G\times F\to G'\times
F'$ is a fully faithful morphism of groupoids.

\bigskip

\section{Deformation functors associated to a DGLA}
\label{sec.3}

Let $L=\oplus L^i$ be a DGLA over a field $\K$ of characteristic 0, we
can define the following three functors:\begin{enumerate}
\item The Gauge functor.
$G_{L}\colon Art_{\K}\to Group$, defined by
$G_{L}(A)=exp(L^0\otimes m_A)$. It is immediate to see that
$G_{L}$ is smooth and homogeneous.

\item The Maurer-Cartan functor. $MC_{L}\colon Art_{\K}\to Set$ defined
by
\[MC_L(A)=\{x\in L^1\otimes m_A\mid dx+\frac{1}{2}[x,x]=0\}.\]
$MC_L$ is a homogeneous functor. If  $L$ is abelian then $MC_L$ is smooth.\\

\item Since $L\otimes m_A$ is a DGLA with  $(L\otimes m_A)^0$
nilpotent, by the result of Section \ref{sec.1} there exists an action of the group
functor $G_{L}$ over $MC_{L}$, we call $Def_{L}=MC_L/G_{L}$ the
corresponding quotient. $Def_L$ is called the {\bf deformation functor}
associated to the DGLA $L$.\\
In general $Def_L$ is not homogeneous but it is only a deformation
functor.
\end{enumerate}

It is also evident that every morphism $\alpha\colon L\to N$ of
DGLA induces morphisms of functors $G_{L}\to G_{N}$,
$MC_{L}\to MC_{N}$. These morphisms are compatible with the gauge action
and therefore induce a morphism between the deformation functors
$Def_L\to Def_M$.\\

We are now ready to compute tangent and complete obstructions for the
above functors.\\
a) $G_{L}$ is smooth, its tangent space is
$G_{L}(\K[\epsilon])=L^0\otimes \K\epsilon$, where $\epsilon^2=0$.\\
b) The tangent space of $MC_L$ is
\[ t_{MC_L}=\{x\in L^1\otimes\K\epsilon \mid dx+\frac{1}{2}[x,x]=0\}
=Z^1(L)\otimes\K\epsilon.\]
If $b\in L^0\otimes\K\epsilon$ and $x\in Z^1\otimes\K\epsilon$ we have
$exp(b)(x+d)=exp(ad(b))(x+d)=x+d+[b,x+d]_d=x+d+db$ and therefore the action
$t_{G_{L}}\times t_{MC_{L}}\to t_{MC_{L}}$ it is given by
$(b,x)\mapsto x+db$.\\
Let's consider a small extension in $Art_{\K}$
\[ e\colon\qquad  0\mapor{}J\mapor{} A\mapor{}B\mapor{}0\]%
and let $x\in MC_L(B)=\{x\in L^1\otimes m_B\mid dx+\dfrac{1}{2}[x,x]=0\}$;
we define an obstruction  $v_e(x)\in H^2(L\otimes J)=
H^2(L)\otimes J$ in
the following way:\\
First take a lifting $\tilde{x}\in L^1\otimes m_A$ of $x$ and consider
$h=d\tilde{x}+\dfrac{1}{2}[\tilde{x},\tilde{x}]\in L^2\otimes J$; we have
\[ dh=d^2\tilde{x}+[d\tilde{x},\tilde{x}]=
[h,\tilde{x}]-\frac{1}{2}[[\tilde{x},\tilde{x}],\tilde{x}].\]
Since $[L^2\otimes J,L^1\otimes m_A]=0$ we have $[h,\tilde{x}]=0$, by
Jacobi identity $[[\tilde{x},\tilde{x}],\tilde{x}]=0$ and then $dh=0$.\\
We define $v_e(x)$ as the class of $h$ in $H^2(L\otimes J)=H^2(L)\otimes
J$; the first thing to prove is that $v_e(x)$ is independent from the
choice of the lifting $\tilde{x}$; every other lifting is of the form
$y=\tilde{x}+z$, $z\in L^1\otimes J$ and then
\[ d\tilde{y}+\dfrac{1}{2}[y,y]=h+dz.\]
It is evident from the above computation that:\begin{enumerate}

\item $(H^2,v_e)$ is a complete obstruction theory for the functor
$MC_L$.

\item If $\phi\colon L\to M$ is a morphism of DGLA then the linear map
$H^2(\phi)\colon H^2(L)\to H^2(M)$ is a morphism of obstruction spaces
compatible with the morphism $\phi\colon MC_L\to MC_N$.
\end{enumerate}

Let's compute the primary obstruction map
$v_{2}\colon Z^1(L)=t_{MC_L}\to H^2(L)$ associated to the
small extension
\[ 0\mapor{}\K\mapor{t^2}\K[t]/(t^3)\mapor{}\K[t]/(t^2)\mapor{}0.\]
Let $x\in MC_L(\K[t]/(t^2))=Z^1(L)\otimes\K t$ and let $\tilde{x}$ be a
lifting of $x$ to $L^1\otimes (t)/(t^3)$, we can choose $\tilde{x}\in
Z^1(L)\otimes (t)/(t^3)$ and therefore the primary obstruction
$v_{2}(x)$ is the class of
$h=d\tilde{x}+\dfrac{1}{2}[\tilde{x},\tilde{x}]=\dfrac{1}{2}[{x},{x}]$
into $H^2(L)\otimes \K t^2$.\\

{\sc Exercise:} Assume either $H^2(L)=0$ or $[L^1,L^1]\subset B^2$,
then $MC_L$ is smooth. If $MC_L$ is smooth then $[Z^1,Z^1]\subset B^2$.\\

{\sc Exercise:} If $[Z^1,Z^1]=0$ then $MC_L$ is smooth.\\

We shall prove later that if $L$ is formal then $MC_L$ is smooth if and
only if $[Z^1,Z^1]\subset B^2$.\\
\medskip\\
c) The tangent space of $Def_L$ is simply the quotient of $Z^1(L)$ under the
translation action of $L^0$ defined by the map $d\colon L^0\to Z^1(L)$ and
then $t_{Def_L}=H^1(L)$.\\
The projection $MC_L\to Def_L$ is smooth and
then induces an isomorphism between universal obstruction theories; by
\ref{2.21} we can define naturally a complete obstruction theory $(H^2(L),
o_e)$ by setting $o_e(x)=v_e(x')$ for every small extension $e$ as above,
$x\in Def_L(B)$ and $x'\in MC_L(B)$ a lifting of $x$.\\
In particular the primary obstruction map $o_2\colon H^1(L)\to H^2(L)$ is
equal to $o_2(x)=\dfrac{1}{2}[{x},{x}]$.\\

\begin{theorem}\label{3.1}
Let $\phi\colon L\to M$ be a morphism of DGLA and denote
by $H^i(\phi)\colon H^i(L)\to H^i(N)$ the induced maps in cohomology.
\begin{enumerate}
\item If $H^1(\phi)$ is bijective and $H^2(\phi)$ injective then the
morphism $Def_L\to Def_N$ is \'etale.\\
\item If, in addition to 1),
the map $H^0(L)\to H^0(N)$ is surjective then
$Def_L\to Def_N$ is an isomorphism.
\end{enumerate}
\end{theorem}

\begin{proof} In case 1) the morphism $Def_L\to Def_M$ is bijective on tangent
spaces and injective on obstruction spaces, by the standard smoothness
criterion it is \'etale.\\
In case 2), since \'etale morphisms are surjective, it is sufficient to
prove that, for every $S\in Art_{\K}$ the map $Def_L(S)\to Def_N(S)$ is
injective.\\
Let $a,b\in MC_L(S)$, as in Section \ref{sec.2} we define the functor
\[ Iso(a,b)\colon Art_S\to Set,\qquad Iso(a,b)(S\mapor{\eta}A)=\{g\in
exp(L\otimes m_A)\mid g*\eta(a)=\eta(b)\}.\]
Since $MC_L$ is homogeneous, $Iso(a,b)$ is a deformation functor with
tangent space $Z^0(L)$ and complete obstruction space $H^1(L)$.\\
Let $K_a\colon Art_S\to Grp$ the group functor defined by
\[ K_a(S\mapor{\eta}A)=
exp(\{[\eta(a),b]+db\mid b\in L^{-1}\otimes m_A\})\subset exp(L^0\otimes 
m_A).\]

The above definition makes sense since it is easy to see that
$\{[\eta(a),b]+db\mid b\in L^{-1}\otimes m_A\}$ it is a Lie subalgebra
of $L^0\otimes m_A$. Since $L^{-1}\otimes m_A\to L^{-1}\otimes m_B$
is surjective
for every surjective $A\to B$, the functor $K_a$ is smooth.\\
Moreover $K_a$ is a subfunctor of $Iso(a,a)$ and therefore acts by right
multiplication on $Iso(a,b)$.
Again by  the result of Section \ref{sec.2} the
quotient functor $Iso(a,b)/K_a$ is a deformation functor with tangent
space $H^0(L)$ and complete obstruction space $H^1(L)$.\\
Now the proof follows as in \ref{2.23};
the morphism
$Iso(a,b)/K_a\mapor{\phi}Iso(\phi(a),\phi(b))/K_{\phi(a)}$ is smooth and
in particular surjective.\end{proof}

\begin{cor}\label{3.2} Let $L\to N$ be a quasiisomorphism of DGLA.
Then the induced morphism $Def_L\to Def_N$ is an isomorphism.\end{cor}
\begin{proof} Evident.\end{proof}

{\sc Exercise:} Let $L$ be a formal DGLA, then $Def_L$ is smooth if and
only if the quadratic map $[,]\colon H^1\times H^1\to H^2$ is zero.\\

\begin{example}\label{3.3}
Let $N^1\oplus B^1= L^1$ be a vector space decomposition
and consider the DGLA  $N=\oplus N^i$ where
\[ \left\{\begin{array}{ll}
N^i=0& \text{ for } i\le 0\\
N^1=N^1&\\
N^i=L^i& \text{ for } i\ge 2\end{array}\right.\]
with bracket and differential induced from $L$. Then the natural inclusion
$N\to L$ gives isomorphisms $H^i(N)\to H^i(L)$ for every $i\ge 1$. In
particular the morphism $Def_N\to Def_L$ is \'etale; note that $Def_N=MC_N$
is homogeneous.\end{example}

{\sc Exercise:} Let $z\in MC_L(S)$, $S\in Art_{\K}$ and consider the
functor $P_z\colon Art_S\to Set$ given by
\[ P_z(A)=\{x\in L^1\otimes m_A\mid dx+[z,x]=0\}.\]
Compute tangent and obstruction spaces of $P_z$.\\

\begin{remark} The abstract
$T^1$-lifting theorem (cf. \cite{FM2}) implies that if the functor
$P_z$ is smooth (informally this means that the linear operator
$d+ad(z)$ has constant rank)
for every $n>0$ and every
$z\in MC_L(\K[t]/(t^n))$ then $MC_L$ is also smooth.
\end{remark}

\begin{cor}\label{3.4} Let $L$ be a DGLA; if $H^0(L)=0$ then $Def_L$ is
homogeneous.\end{cor}

\begin{proof} It is sufficient to apply \ref{3.1} at the Example \ref{3.3}.\end{proof}

\begin{example}\label{3.5} {(The simplest example of deformation
problem governed by a DGLA)}\\
Let $V$ be a $\K$-vector space and $J\colon V\to V$ linear
such that $J^2=-I$ (if $\K=\R$ then $J$ is a complex structure on $V$).
Let $M=\Hom(V,V)$ and $d,\delta\colon M\to M$ defined by
\[ d(A)=JA+AJ,\qquad \delta(A)=JA-AJ,\qquad A\in M.\]
Since $char\K\not=2$ it is easy to prove that 
$\ker d=\operatorname{Image}(\delta)$,
$\operatorname{Image}(d)=\ker\delta$.\\
Consider the DGLA $L=\oplus_{n\ge 1}M[-n]$ with differentials
\[ 0\mapor{}M[-1]\mapor{d}M[-2]\mapor{\delta}M[-3]
\mapor{d}M[-4]\mapor{\delta}\ldots\]
and the bracket defined by
\[ [A,B]=AB-(-1)^{\overline{A}\overline{B}}BA.\]
Clearly $H^i(L)=0$ for every $i\not=1$ and then the associated functor
$Def_L=MC_L$ is
smooth homogeneous with tangent space $\ker{d}$.\\
The elements of $MC_L(A)$ are in natural bijection with the deformation
of the solution of $J^2=-I$
over the Artinian ring  $A$; in fact given $H\in M\otimes m_A$ we have
$H\in MC_L(A)$ if and only if $(J+H)^2=-I$.
Since the gauge group is trivial, this example makes sense
over arbitrary fields of $char\K\not=2$.\end{example}

\bigskip

\section{The Kuranishi map and the Kuranishi functor}

By a well known theorem of Schlessinger \cite{Sch}, if $F\colon Art_{\K}\to Set$
is a deformation functor
with finite dimensional tangent vector space then there exists a prorepresentable
functor $h_R$ and an \'etale morphism $h_R\to F$.\\
In particular Schlessinger theorem applies to the functor $Def_L$ for
every DGLA $L$ such that $H^1(L)$ is finite dimensional. In this section
we give a explicit construction of a prorepresentable functor $h_R$ and
of an \'etale map $h_R\to Def_L$ by introducing Kuranishi maps.\\

Let $L=\oplus L^i$ be a DGLA, for every $i\in \Z$ choose direct sum
decomposition
\[ Z^i=B^i\oplus \sH^i,\qquad L^i=Z^i\oplus C^i.\]
Let $\delta\colon L^{i+1}\to L^i$ be the linear map composition of
\begin{enumerate}
\item The projection $L^{i+1}\to B^{i+1}$ of kernel $C^{i+1}\oplus
\sH^{i+1}$.
\item $d^{-1}\colon B^{i+1}\to C^i$,  the inverse of $d$.
\item The inclusion $C^i\to L^i$.
\end{enumerate}

We note that $x\in \sH^i$ if and only if $dx=\delta x=0$ and
$d\delta+\delta d=Id-H$, where $H\colon L^1\to \sH^i$ is the projection of
kernel $B^i\oplus C^i$.\\

Using the same notation of \ref{2.5}.1,
for a vector space $V$ we call $\widehat{V}\colon Art_{\K}\to Set$ the
homogeneous functor $\widehat{V}(A)=V\otimes m_A$.\\

\begin{definition}\label{4.1}
The Kuranishi map $F\colon \widehat{L^1}\to \widehat{L^1}
$ is the morphism of functors given by
\[ x\in L^1\otimes m_A,\qquad F(x)=
x+\dfrac{1}{2}\delta[x,x]\in L^1\otimes m_A.\]
\end{definition}

\begin{lem}\label{4.2} The Kuranishi map $F$ is an isomorphism of
functors.\end{lem}

\begin{proof} $F$ is clearly a morphism of functors, since $\widehat{L^1}$ is
homogeneous, by \ref{2.12} it is sufficient to prove that $F$ is \'etale.\\
By construction $F$ is the identity on the tangent space of $\widehat{L^1}$.
Since $\widehat{L^1}$ is smooth, the trivial obstruction theory $(0,v_e)$ is
complete; by standard smoothness criterion \ref{2.18} $F$ is
\'etale.\end{proof}

\begin{definition}\label{4.3}
The Kuranishi functor $Kur\colon Art_{\K}\to Set$ is
defined by
\[ Kur(A)=\{ x\in \sH^1\otimes m_A\mid H([F^{-1}(x),F^{-1}(x)])=0\}.\]
\end{definition}

In other words $Kur$ is the kernel the morphism of homogeneous functors
$\widehat{q}\colon \widehat{\sH^1}\to\widehat{\sH^2}$
induced by the map of formal pointed schemes $q\colon \sH^1\to \sH^2$ which is
composition of: the inclusion
${\sH^1}\to {L^1}$, the inverse of $F$, the bracket $[,]\colon
{L^1}\to {L^2}$ and the projection ${L^2}\to {\sH^2}$.

\begin{lem}\label{4.4} $Kur$ is homogeneous; if moreover the dimension of
$H^1(L)$ is finite then $Kur=h_R$ is prorepresented by a $\K$-algebra
$R\in \widehat{Art}_{\K}$ such that $Spec(R)=(q^{-1}(0),0)$ as formal
scheme.\end{lem}

\begin{proof} Exercise (easy).\end{proof}

\begin{lem}\label{4.5} Let $x\in \widehat{L^1}(A)$, $y\in \widehat{L^i}(A)$, $c\in \K$
such that
$y=c\delta[y,x]$, then $y=0$.\end{lem}

\begin{proof} Exercise (easy).\end{proof}

Clearly $MC_L$ and $Kur$ can be considered as subfunctor of $\widehat{L^1}$.

\begin{prop}\label{4.6} The isomorphism $F\colon \widehat{L^1}\to \widehat{L^1}$
induces an isomorphism
\[ F\colon MC_L\cap \widehat{(C^1\oplus\sH^1)}\to Kur.\]
\end{prop}

\begin{proof} If $x\in MC_L(A)\cap (C^1\oplus\sH^1)\otimes m_A$,
then $dx+\dfrac{1}{2}[x,x]=0$, $\delta x=0$ and $H([x,x])=H(-2dx)=0$.
Moreover we have:\\
\[ [x,x]=[x,x]-H([x,x])=d\delta[x,x]+\delta d[x,x]=d\delta[x,x]\]
\[ \delta F(x)=\delta x+\dfrac{1}{2}\delta^2[x,x]=\delta x\]
\[ dF(x)=dx+\dfrac{1}{2}d\delta[x,x]=dx+\dfrac{1}{2}[x,x]=0.\]
This prove that $F(x)\in Kur(A)$.\\
Conversely if $x\in L^1\otimes m_A$ and $F(x)\in Kur(A)$ then
$\delta x=\delta F(x)=0$, $H([x,x])=0$.\\
$dx+\dfrac{1}{2}[x,x]=dx+\dfrac{1}{2}(d\delta+\delta d)[x,x]=
dF(x)+\dfrac{1}{2}\delta d[x,x]=\dfrac{1}{2}\delta d[x,x]$.\\
Therefore it is sufficient to prove that $\delta d[x,x]=0$;
\[ \delta d[x,x]=2\delta[dx,x]=-\delta [d\delta[x,x],x]=
-\delta[[x,x]-H([x,x])-\delta d[x,x],x].\]
Since $H([x,x])=0$ and by Jacoby $[[x,x],x]=0$ we have
\[ \delta d[x,x]=\delta[\delta d[x,x],x]\]
and the proof follows from Lemma \ref{4.5}.\end{proof}

Proposition \ref{4.6} says that, under the Kuranishi map $F$, the Kuranishi
functor $Kur$ is isomorphic to $MC_N$, where $N=\oplus N^i$ is the DGLA
defined by
\[ \left\{\begin{array}{ll}
N^i=0& \text{ for } i\le 0\\
N^1=C^1\oplus \sH^1&\\
N^i=L^i& \text{ for } i\ge 2\end{array}\right.\]
By Example \ref{3.3} the morphism $MC_N\to Def_L$ is \'etale. In conclusion what
we have proved is the following

\begin{theorem}\label{4.7} For every DGLA $L$ the morphism
$Kur\mapor{F^{-1}} MC_L\to Def_L$ is \'etale.
\end{theorem}

{\sc Exercise:} In the notation of \ref{3.3}, if $H^1(L)$ is finite dimensional
then the isomorphism class of $MC_N$ is independent from the choice of
the complement $N^1$. (This is always true but a proof of this fact
without assumptions on $H^1(L)$ requires a nontrivial application of the
factorization theorem \cite[6.2]{FM1}).

\bigskip

\section{Homotopy equivalence versus gauge equivalence}

We introduce here a new equivalence relation in the set of
solutions of the Maurer-Cartan equation in a differential graded Lie
algebra. Let $(L, [,],\delta)$, $L=\oplus L^i$, be a fixed DGLA and let
$\K[t,dt]=\K[t]\oplus\K[t]dt$ be the differential graded algebra
with $t$ of degree 0, $dt$ of degree 1, endowed with the expected
differential $d(p(t)+q(t)dt)=\dfrac{dp(t)}{dt}dt$.\\
We denote by $\Omega=L\otimes_{\K}\K[t,dt]=\oplus \Omega^i$; it is a
DGLA with $\Omega^i=L^i[t]\oplus L^{i-1}[t]dt$, differential
$\delta_\Omega(a(t)+b(t)dt)=\delta a(t)+(-1)^{\overline{a(t)}}
\dfrac{\partial a(t)}{\partial t}dt
+\delta b(t)dt$ and bracket
\[ [a(t)+b(t)dt,p(t)+q(t)dt]=
[a(t),b(t)]+[a(t),q(t)]dt+(-1)^{\overline{p(t)}}[b(t),q(t)]dt.\]
For every $s\in\K$ there exists an evaluation morphism of DGLA
$v_s\colon \Omega\to L$, $v_s(t)=s$, $v_s(dt)=0$. In particular we have
morphisms of functors of Artin rings $v_s\colon MC_\Omega\to MC_L$ which
are left inverses of the natural inclusion $MC_L\subset MC_\Omega$.

\begin{lem}\label{5.1} The morphism $v_s\colon MC_\Omega\to MC_L$ is smooth
for every $s\in\K$.\end{lem}

\begin{proof} We first give a more explicit description of the solutions of the
Maurer-Cartan equation in $\Omega$. Take $\omega\in\Omega^1$,
$\omega=a(t)+b(t)dt$, we have
\[ \delta_\Omega\omega+\dfrac{1}{2}[\omega,\omega]=
\delta a(t)+\dfrac{1}{2}[a(t),a(t)]+
\left(-\dfrac{\partial a(t)}{\partial t}+\delta
b(t)+[a(t),b(t)]\right)dt\]
and then $\omega\in MC_\Omega$ if and only if
\begin{enumerate}
\item $a(t)\in MC_L$ for every $t\in\K$, and
\item $\dfrac{\partial a(t)}{\partial t}=\delta b(t)+[a(t),b(t)]$.
\end{enumerate}

Let $0\mapor{}J\mapor{}A\mapor{\pi}B\mapor{}0$ be a small extension of
Artin rings, $\omega_B=a_B(t)+b_B(t)dt\in MC_\Omega(B)$
and $a_A(s)\in MC_L$ such that $\pi(a_A(s))=a_B(s)$.\\
Let $\tilde{a}_A(t)$, $b_A(t)$ be liftings of $a_B(t)$, $b_B(t)$ in
$L[t]\otimes m_A$ such that $\tilde{a}_A(s)=a_A(s)$; since $Jm_A=0$ the
bracket $[\tilde{a}_A(t),b_A(t)]$ is independent from the choice of the
liftings.\\
Setting
\[ \gamma(t)=\delta b_A(t)+[\tilde{a}_A(t),b_A(t)]-
\dfrac{\partial \tilde{a}_A(t)}{\partial t}
\in L^1[t]\otimes J\]
and
\[ a_A(t)=\tilde{a}_A(t)+\int_{s}^t\gamma(\tau)d\tau,\]
it is evident that $\omega_A=a_A(t)+b_A(t)dt$ satisfies the differential
equation $ii)$ and $v_s(\omega_A)=a_A(s)$. Moreover
$\delta a_A(t)+\dfrac{1}{2}[a_A(t),a_A(t)]\in L^2[t]\otimes J$ and then
\begin{multline*} 
\dfrac{d~}{dt}(\delta a_A(t)+\dfrac{1}{2}[a_A(t),a_A(t)])=
\delta \dfrac{d a_A(t)}{dt}+\left[\dfrac{d a_A(t)}{dt}, 
a_A(t)\right]\\
=\delta [a_A(t), b_A(t)]+[\delta b_A(t)+[a_A(t),b_A(t)],a_A(t)]\\
=\left[-\dfrac{1}{2}[a_A(t),a_A(t)], b_A(t)\right]+
[[a_A(t), b_A(t)], a_A(t)]=0.\end{multline*}
By assumption $a_A(s)\in MC_L$ and then $a_A(t)\in MC_L$ for every
$t\in\K$.\end{proof}

\begin{definition}\label{5.2}
We shall say that $x\in MC_L(A)$ is homotopic to $y\in MC_L(A)$ if
there exists $\omega\in MC_\Omega$ such that $x=v_0(\omega)$,
$y=v_1(\omega)$; it is easy to see that the relation ``$x$ is homotopic
to $y$'' is reflexive and symmetric (cf. exercise below). The equivalence
relation generated is called {\it homotopy equivalence}.\end{definition}

{\sc Exercise:} There exists a natural action of the group of affine
isomorphisms of $\K$ on $MC_\Omega$.\smallskip\\
Let's momentarily denote by $F_L\colon Art_{\K}\to Set$ the quotient
functor of $MC_L$ by the homotopy equivalence (at the end of the trip
$F_L$ will be equal to $Def_L$).\\

\begin{prop}\label{5.3} $F_L$ is a deformation functor and the natural
projection $MC_L\to F_L$ is smooth.\end{prop}

\begin{proof} Easy exercise (Hint: use that $v_0,v_1$ are smooth left inverses
of the natural inclusion $MC_L\subset MC_\Omega$.)\end{proof}

\begin{lem}\label{5.4} The tangent space of $F_L$ is the quotient of
$t_{MC_L}$ by the image of the linear map $v_1-v_0\colon t_{MC_\Omega}\to
t_{MC_L}$.\end{lem}

\begin{proof} It is evident that if $x\in t_{MC_L}$ is homotopy equivalent to
$y$
then $x-y$ belongs to the image of $v_1-v_0$. Conversely let $x\in
t_{MC_L}$, $\omega\in t_{MC_\Omega}$, $h=v_1(\omega)-v_0(\omega)$ and
$y=x+y$.\\
We can write $x=v_0(x-v_0(\omega)+\omega)$, $y=v_1(x-v_0(\omega)+\omega)$
and therefore $x$ is homotopy equivalent to $y$.\end{proof}

We note that $\omega=a(t)+b(t)dt\in t_{MC_\Omega}$ if and only if $\delta
a(t)=0$ and $\dfrac{d a(t)}{dt}=\delta b(t)$; in particular
\[v_1(\omega)-v_0(\omega)=a(1)-a(0)=\int_0^1 \delta b(t)dt\in B^1(L).\]
On the other hand for every $b\in L^0$, $\delta(b)t+bdt\in t_{MC_\Omega}$ and
therefore the image of $v_1-v_0$ is exactly $B^1(L)$.\\

\begin{theorem}\label{5.5} $x\in MC_L(A)$ is homotopic to $y\in MC_L(A)$ if and
only if $x$ is gauge equivalent to $y$.\end{theorem}

Before proving Theorem \ref{5.5} we need some preliminary results about the exponential
map: let $N$ be a nilpotent Lie algebra and let $a\to e^a$ be its
exponential map $N\to exp(N)$.\\
Given a polynomial $p(t)\in N[t]$ we have
\[ e^{p(t+h)}e^{-p(t)}=e^{h(p'(t)+\gamma_p(t))+h^2\eta_p(t,h)}\]
where $p'(t)=\dfrac{\partial p(t)}{\partial t}$ and $\gamma_p,\eta_p$ are
polynomials determined exactly by Taylor expansion
$p(t+h)=p(t)+hp'(t)+h^2\cdots$ and Campbell-Baker-Hausdorff formula. For
later use we point out that
$\gamma_p(t)$ is a linear combination with rational coefficients of terms
of the following types:
\[ ad(p(t))^{\alpha_1}ad(p'(t))^{\beta_1}\cdots
ad(p(t))^{\alpha_n}ad(p'(t))^{\beta_n}Q\]
with $Q=p(t),p'(t)$ and $\sum \alpha_i+\sum\beta_i>0$.

\begin{lem}\label{5.6} In the notation above, for every $b(t)\in N[t]$ there
exists an unique polynomial $p(t)\in N[t]$ such that $p(0)=0$ and
$p'(t)+\gamma_p(t)=b(t)$.\end{lem}

\begin{proof} Let $N=N^{(1)}\supset N^{(2)}\supset\cdots\supset N^{(m)}=0$ be the
descending central series, $N^{(i+1)}=[N,N^{(i)}]$, and denote
$N_i=N/N^{(i+1)}$. By construction $N_{i+1}\to N_i$ is a central
extension for every $i\ge 0$.\\
The lemma is trivially true over $N_0$; assume that $p_i(t)\in N_i(t)$ is
the unique solution of the differential equation
$p_i'(t)+\gamma_i(t)=b(t)$ ($\gamma_i=\gamma_{p_i}$) on the vector
space $N_i[t]$.\\
Let $\tilde{p}_{i+1}(t)\in N_{i+1}$ be a lifting of $p_i(t)$ such that
$\tilde{p}_{i+1}(0)=0$, as the kernel $K_{i+1}$ of the projection $N_{i+1}\to
N_i$ is contained in the centre the polynomial
$\gamma_{i+1}=\gamma_{\tilde{p}_{i+1}}$ does not depend from the choice
of the lifting.\\
We introduce the polynomials
\[ \chi(t)=b(t)-\tilde{p}_{i+1}'(t)-\gamma_{i+1}(t)\in K_{i+1}[t]\]
\[ {p_{i+1}}(t)=\tilde{p}_{i+1}(t)+\int_0^t\chi(\tau)d\tau\in 
N_{i+1}[t]\]
It is clear that $p_{i+1}(t)$ is the unique solution of the differential
equation.\end{proof}
If $N\times V\mapor{[,]}V$ is a representation of $N$ into a vector space
$V$, for every $v\in V$ we have
\[ \dfrac{d~}{dt}(e^{p(t)}v)=[p'(t)+\gamma_p(t),e^{p(t)}v]\]
As an immediate consequence of 5.6 we have

\begin{prop}\label{5.7} Let $N$ be a nilpotent Lie algebra, for every $b(t)\in
N[t]$ there exists unique $p(t)\in N[t]$ such that $p(0)=0$ and
for every representation $N\times V\mapor{[,]}V$, for every $v\in V$ holds
\[ \dfrac{d~}{dt}(e^{p(t)}v)+[b(t),e^{p(t)}v]=0.\]
\end{prop}

\begin{proof} Evident.\end{proof}

\begin{proof}[Proof of 5.5] The gauge equivalence and the notion of homotopy
extends naturally to the space $V=(L^1\otimes
m_A)\oplus\K\delta$:\begin{enumerate}

\item $x,y\in V$ are gauge equivalent if and only if there exists
$p\in L^0\otimes m_A$ such that $y=e^px$.

\item $x\in V$ is homotopic to $y\in V$ if and only if there
exists $a(t)\in V[t]$, $b(t)\in L^0\otimes m_A[t]$ such that
$a(0)=x$, $a(1)=y$ and $\dfrac{d a(t)}{dt}+[b(t),a(t)]=0$.
\end{enumerate}
Now \ref{5.5} follows from \ref{5.7}.\end{proof}

\begin{remark} The notion of homotopy equivalence extends naturally to
$L_{\infty}$-algebras see \cite{Can},
while the gauge action requires the
structure of DGLA; this
motivates this section and the redundancy of proofs given here.
\end{remark}

\end{document}